\documentclass[twocolumn,10pt]{articlefederico}

\usepackage{amsthm}
\usepackage{amsmath}
\usepackage{amssymb}
\usepackage{color}

\usepackage{graphicx}

\begin{document}
\title{\textsf{Todxs cuentan in ECCO: community and belonging in mathematics}}
\author{\textsf{Federico Ardila--Mantilla\footnote{\noindent San Francisco State University, Universidad de Los Andes}\,\,\, and Carolina Benedetti\footnote{\noindent Universidad de Los Andes}}}

\date{}
\maketitle

\begin{quote}
\emph{Aprend\'{\i} que para uno encontrarse tiene que buscar en la ra\'{\i}z. [...]
Aprend\'{\i} que no soy s\'olo yo, y que somos muchos m\'as.}

Grupo Bah\'{\i}a
\end{quote}

\section{\textsf{What is ECCO?}}

The Encuentro Colombiano de Combinatoria is a biannual gathering of students and researchers from  Colombia, Latin America, California, and many other places. 
It is a two-week long summer school, featuring mini-courses by experts, collaborative problem workshops, research talks and posters, open problem sessions, a discussion panel, a hike, and visits to some of Colombia's legendary salsa clubs. 
It is also much more than a summer school, and we hope to capture a bit of its spirit in these pages.

ECCO is designed to give \textbf{every} participant opportunities to interact closely with people at all stages of the mathematical career. We do our best to build a very professional and very warm atmosphere. We are collaborators and we are also a community.

The Encuentro started as a small gathering for combinatorics students in Colombia and the San Francisco Bay Area. They had taken classes together, as part of the SFSU-Colombia Combinatorics Initiative described in \cite{TodosCuentan}, and it had become clear that they wanted to meet in person, build closer ties, and find ways to collaborate.

Since then, ECCO has broadened and gained a strong reputation. Students from many different countries now attend, and combinatorics experts also ask us to participate. We communicate our goals clearly. This is not a regular conference; it is a school and an \textbf{encuentro}: a coming together. 
We ask these experts to do problem sets with the students, to present research questions that they would like help with, to offer advice, and to join the dance floor at some point. They have been wonderfully helpful and inspiring mentors, they have recruited students, and -- perhaps most meaningfully to us -- several have mentioned that their experiences at ECCO have influenced their work at their home institutions.

As one becomes more experienced organizing events, one becomes more conscious of their shortcomings. ECCO is certainly an imperfect event. After seventeen years, it is still under construction, and we hope it continues to be.
But ECCO has been tremendously inspiring and energizing to us, and has taught us a lot about what it might mean to truly find community and belonging in a mathematical space. The goal of this article is to share a few of the lessons that we have learned from helping to build it.

\section{\textsf{Community Agreement, Part 1}}

When prospective participants are applying to ECCO, they encounter our Community Agreement. The first part reads:

\begin{quote}
\emph{\textbf{A rewarding experience for all.} The Encuentro Colombiano de Combinatoria aims to offer a rewarding, challenging, supportive, and fun experience to every participant. We will build that rich experience together by devoting our strongest available effort to all ECCO activities. You will be challenged and supported. Please be prepared to take an active, critical, patient, and generous role in your own learning and that of the other participants.}
\end{quote}

When we meet in person, we start ECCO by reminding everyone about this agreement.
We ask people to get in pairs, read it out loud to each other, and spend a few minutes discussing it: What stands out to you about this agreement? What can it look like to put it in practice?

We're not gonna lie. While some participants jump right in, many look confused, and if we are reading their body language correctly, a few seem to think: \emph{I can't believe you are asking me to do this; what am I, a kindergartener?} But we insist. Everyone participates. 

To initiate a dialogue, we ask each group to underline a few words in the agreement that  resonate with them, and share them with everyone. Some are excited that they will be challenged; some that they will be supported; some point out that the combination is crucial. 
We discuss how to be productively critical of each other's work, and what generosity might mean in a mathematical setting. We talk about how sometimes we are very good at being patient with others, but not so good at being patient with ourselves.

We wrote this agreement to communicate, from day one, the kind of space we are trying to build collectively. Johan, one of the participants of  D\'ias de Combinatoria\footnote{The D\'{\i}as summer school is one of the offsprings of ECCO.},  shared with us an experience that became an unforeseen consequence of the agreement. He told us that reading it on the webpage of D\'ias was the push he needed to apply, and to attend; for the first time, he felt he was welcome in an event like this.

\section{\textsf{Community Agreement, Part 2}}

The second part of the community agreement reads:\footnote{This part of the agreement was based on a code of conduct written by Ashe Dryden; she has compiled many valuable resources on this topic. \cite{Dryden}}

\begin{quote}
\emph{\textbf{A welcoming experience for all.}
 ECCO is committed to creating a professional and welcoming environment that benefits from the diversity of experiences of all its participants. We will not tolerate any form of discrimination or harassment. We aim to offer equal opportunity and treatment to every participant regardless of their mathematical experience, gender identity, nationality, race or ethnicity, religion, age, marital status, sexual orientation, disability, or any other factor. \\ 
Behavior or language that is welcome or acceptable to one person may be unwelcome or offensive to another. Consequently, we ask you to use extra care to ensure that your words and actions communicate respect for others. This is especially important for those in positions of authority or power, since individuals with less power have many reasons to fear expressing their objections regarding unwelcome behavior.\\
If a participant engages in discriminatory or harassing behavior, ECCO organizers may take any action they deem appropriate, from warning the offender to immediately expelling them from the event. \\
If you are being harassed, you feel uncomfortable with the way you are being treated, you notice that someone else is being harassed, or you have any other concerns, please contact Carolina
Benedetti or \linebreak Federico Ardila immediately. If you prefer not to speak in person, you may e-mail us (anonymously, if you wish) at the account \underline{\hspace{1.5cm}}@\underline{\hspace{1cm}}.\underline{\hspace{.5cm}}, which only Federico and Carolina have access to.}
\end{quote}

Again, we make sure everyone actively engages with this text, reading it out loud in pairs and discussing it, awkward as they might find that. We are direct: social events are an essential part of ECCO, and we explicitly ask participants not to use them as excuses for romantic advances. We bring together more than 100 strangers from many different cultural backgrounds for an intense shared experience; it is essential to have an agreement that clarifies expectations, and gives the organizers the power to react to potential incidents.

We have co-created these agreements with our students at SFSU, ECCO, and D\'ias de Combinatoria over the last few years. We prefer to call it a \emph{community agreement} instead of a \emph{code of conduct}, and take the time to discuss it, so it does not just feel like an externally imposed code of conduct that they must obey, or worse, bureaucratic fine print.
Our hope is to reach a collective agreement that is actually ours, that everyone is committed to.

In the `Any additional comments?' question on the exit survey of ECCO 2018, 
almost all participants who identified as women and/or LGBTQ+
praised the community agreement, and several said they would like to have one in all math events. Two participants wrote:

\begin{quote}
\emph{I thought the community agreement was an excellent idea. The openness it allowed us to have and to make a giant community out of everyone made the conference very special. I felt I could finally be myself after years of feeling caged in.} 
\bigskip \\
\emph{We made an agreement to acknowledge each other's differences and try our best to create a positive experience for everyone and it worked! We came together and did math without fear or judgment. It was so much fun! I think that the community agreement and the leadership of the organizers, TAs, and Colombians were driving forces behind making that possible. We all played a part by putting our hearts into creating the environment we were longing for. I left feeling fired up about bringing ECCO home with me. I would love it if all of my classes started off with a community agreement at the beginning of the semester.} 
\end{quote}

A senior participant later told us: ``\emph{I was very surprised at first, and looked at [the agreement] as an oddity. Then I remembered what it was like being a grad student at conferences, and all the weird guys I had to avoid. So I figured, yeah, why not?}" 
Another participant, who had been assaulted in a mathematical space before, told us that she simply does not attend conferences that do not have a plan to ensure her safety.

\section{\textsf{Breaking power structures.}}

In any group of people there is a hidden power structure that influences who leads the discussion, who participates, whose voices are listened to, whose ideas are seen as important. 
Our activities are most successful -- for teachers and for students -- when we are able to disrupt those power structures as much as possible, when every participant feels that their presence is important and their thoughts are valuable. We try to do this constantly, in several ways; a particularly successful one occurs outside of the classroom. 

\begin{figure}[h]
\begin{center}
    \includegraphics[width = 8cm]{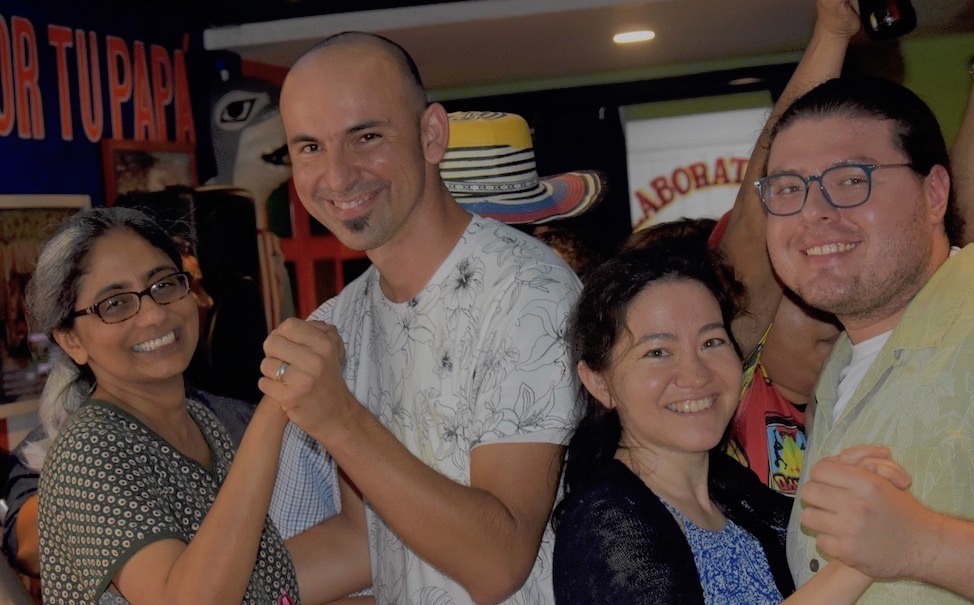}
	\caption{Lecturers and students on the dance floor.}
	\vspace{-.5cm}
\end{center}	
\end{figure}
On Saturday nights, ECCO moves to the dance floor of the best salsa club we can find. The truth is that many of our international experts look a bit intimidated when they first walk in. For most of them, this is not the kind of place they visit often, if ever. Few people at the discoteca look like them; they might feel like they don't really belong there. 

Very soon, the students approach them and invite them to dance. They don't accept ``I don't know how to dance" for an answer; they teach them, patiently, kindly, from the beginning, or just persuade them to dance as they will.

We won't pretend our guests become expert dancers overnight; that really does not matter. But  they always seem really grateful to the students who make sure they are comfortable, who guide them through a few steps, who probably help them find a bit of freedom inside their body. Some of us have known these professors for years, and we get to see a smile that they have never shown us before.

We like to ask our course instructors to keep in mind the feeling of discomfort they might have had entering the discoteca, and the feeling of growth and joy they hopefully had walking out. Many ECCO students -- who have never met so many accomplished mathematicians, who may have never attended a math conference before --  are probably feeling a similar discomfort  when they walk into the classroom; we want them to have that sense of belonging, growth, and empowerment when they leave. Since the professor was vulnerable in front of the student, the student can more comfortably say ``I don't understand, can you explain this to me?" when needed. Since the student showed generosity and patience on the dance floor, the professor naturally shows a similar generosity and patience in the classroom. Dancing serves as a way of illustrating how classrooms can be made more joyful and empowering spaces to everyone.

\begin{quote}
\emph{We must return ourselves to a state of embodiment in order to deconstruct the way power has been traditionally orchestrated in the classroom. }

bell hooks \cite{hooks}
\end{quote}

The dance floor is one of the most democratic spaces of the tremendously unequal societies we live in. At ECCO it is a place of joy, and also a place of pedagogy, for professors and students alike.

\section{\textsf{Problem workshops: thinking simply about deep things}}

Mathematically, ECCO aims for a \textbf{low floor, high ceiling} approach. We want the courses and activities to be designed so that everyone is able to engage with them at some level, and no one runs out of questions to explore. \textbf{Every} participant should find interesting things to learn. This is perhaps best exemplified in the way that problem workshops are structured. 

Each minicourse meets four times, and each 
60-minute class meeting is followed by a 90-minute problem workshop. People self-identify their level of expertise, and we split them into groups as heterogeneously as we can. A typical group will include a professor or postdoc, a graduate student in combinatorics, and two undergraduates with scarce combinatorial experience. Many participants speak very little English, and many speak very little Spanish, so everyone has something to learn and something to teach. We offer materials in Spanish or English, and the unofficial language of mathematical discussions is Spanglish. People are welcome to use the language they wish; interestingly, many choose to communicate in a foreign language for the first time; this is their opportunity to try it.

\begin{figure}[h]
\begin{center}
    \includegraphics[width = 8cm]{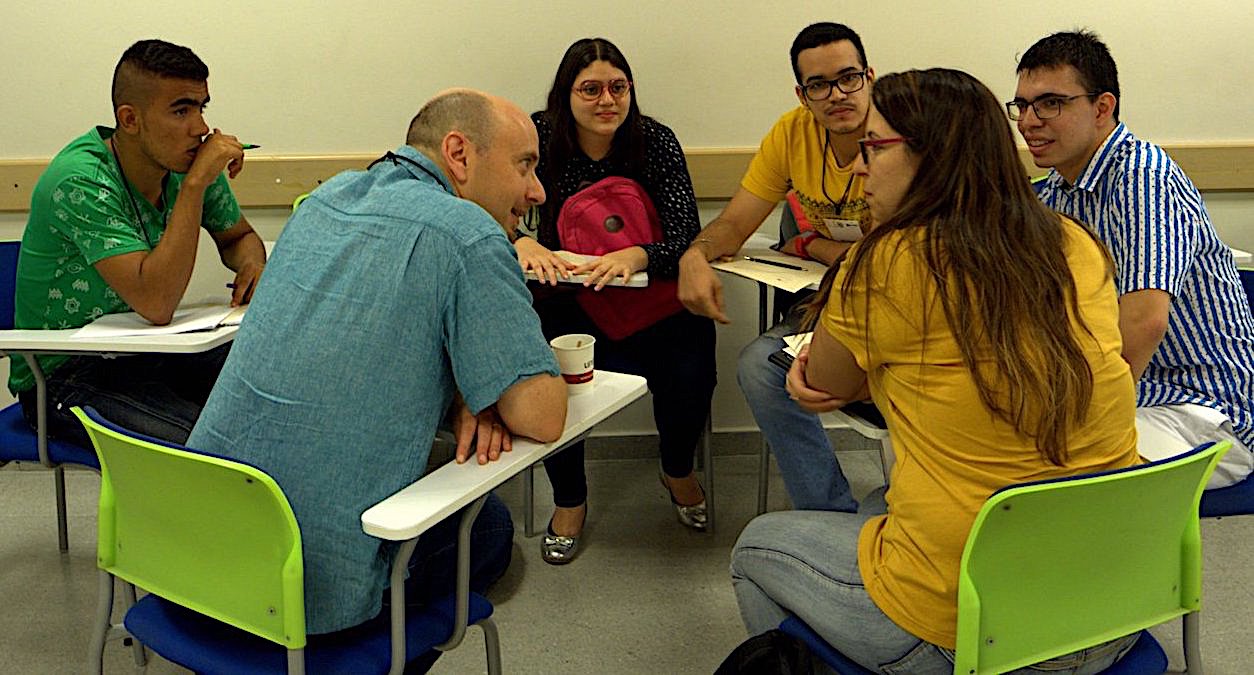}
	\caption{Collaboration in a problem workshop.}
		\vspace{-.5cm}
\end{center}	
\end{figure}

The first problems on each list ask people to carry out a small example, to ensure that everyone understands the key constructions or results in the class. The last few problems on the list can be very challenging, and may take  days or weeks to solve. 
We ask each group to keep in mind our community agreement: how can they make the problem session rewarding, challenging, supportive, and fun for every participant? The result has always exceeded our expectations.

We realize this approach is unusual. Occasionally, it faces some resistance. A few of the more experienced participants have asked: ``Why don't you let the beginners work on the easy problems together, and we can focus on the hardest problems?" But this is how these experts have been operating for most of their career; why not learn something new?

It is very rare for an undergraduate to get to collaborate with an expert of one field on questions about a different field, and see:  \emph{Experts struggle too! How do they productively struggle?} These are very valuable lessons for the undergraduates. It is also very rare for an expert to get to collaborate with a relative newcomer to mathematics, as equals. When they find a way to do it, they inevitably deepen their understanding of the subject. 

The last few minutes of the problem workshop are spent sharing solutions. We ask the people who are usually very comfortable speaking up to make space for others.
We invite the least experienced or the least vocal participants to present their work; they are the ones who can grow the most from doing so, and with the right atmosphere and maybe a bit of extra encouragement, they are usually  happy to speak.

Andr\'es Vindas--Mel\'endez, who was a Master's student when he first attended ECCO, described his experience  \cite{Vindas}:

\begin{quote}
\emph{The exercises were mathematically meaningful, but what is noteworthy is that all group members played an active role in reaching a solution and understanding of the concepts. I observed that the more experienced mathematicians went directly to thinking about the abstraction of the problems, where the younger students emphasized a more concrete approach to exemplify the theory occurring in the problem. Of course both ways of thinking are valuable.}
\end{quote} 

This reminds us of Gelfand's request when encountering a new mathematical topic: 

\begin{quote}
\emph{Explain this to me in a simple example; the difficult example I will be able to do on my own.} 

Israel M. Gelfand \cite{Gelfand}
\end{quote}

\noindent Satisfying this request can be very challenging for beginners and experts alike, and it can also be surprisingly rewarding and enlightening. The celebrated Ross Mathematics Program extols the value of \textbf{thinking deeply about simple things}. We agree wholeheartedly, and propose the counterpart as well: there is a tremendous amount to be learned from \textbf{thinking simply about deep things}.

\begin{figure}[h]
\begin{center}
    \includegraphics[width = 8cm]{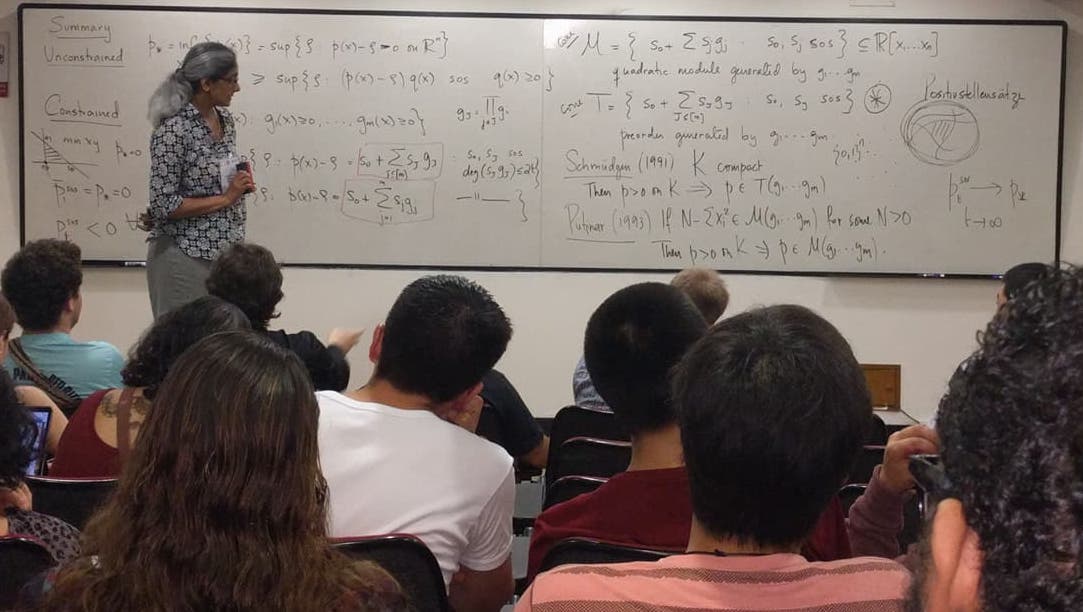}
	\caption{Rekha Thomas teaches optimization.}
		\vspace{-.5cm}
\end{center}	
\end{figure}

\section{\textsf{Universal design.}}

Mathematicians from overrepresented groups in mathematics often ask ``Why do you need these math conferences for minorities? Don't we all do the same mathematics?" To try to answer, allow us to digress for a moment.

Let us share an embarrassing confession: The first time we got in a public bus in North America, and someone in a wheelchair got on, we could not believe our eyes. Are all 50 of us really going to wait all this time for one person to get on? Did the city really spend all this money putting all this equipment on every bus for such a small percentage of the population? 
We both grew up riding the ramschackle buses of Bogot\'a, jumping in and out of them while they were still in motion, collecting frequent minor bruises along the way. We should have known better.

The term ``universal design,"\footnote{We would like to thank May-Li Khoe for teaching us about universal design and its wide applicability.} coined by architect Ron Mace, describes the concept of designing all products and built environments to be aesthetic and usable to the greatest extent possible by everyone, regardless of their age, ability, or status in life. What may be unintuitive about universal design is that, what may seem like designing for a small minority, ends up being a better design for the majority. In fact, once it becomes widely used, it is no longer seen as serving ``special" needs. 

We often forget that sidewalk ramps were installed in every US city thanks to the Americans with Disabilities Act of 1990, after decades of activism. They were originally designed for people who use wheelchairs to go on and off sidewalks easily. Today, everyone uses them: a kid on a tricycle, a parent with a stroller, a traveller with a suitcase, a skate boarder, or the two of us when we were dealing with injuries. Everyone benefits from them.

With this in mind, let us propose an analogy:

\begin{quote}
\emph{Mathematics education cannot truly improve until it adequately addresses the very students who the system has most failed. [...] We need a central focus on students who are Latinx, Black, and Indigenous [...], developing practices and measures that feel humane to those specific communities as a means to guide the field.}

Rochelle Guti\'errez \cite{Gutierrez}
\end{quote}

\noindent
This does not come naturally, or without some opposition. 
Sexism, racism, classism, and centralism often lead to a small, homogeneous group of students being valued more than the rest, tacitly or explicitly. Academic elitism centers the voices and interests of the ``top" students from the ``top" schools, whatever ``top" might mean. Furthermore, in Colombia, we always seem to put the needs of our foreign guests above our own. 

At ECCO we are intentional and unapologetic about focusing on the needs and interests of the local students, the less experienced students, and the students from regional universities that have less access to activities like this. It is our belief, and our experience, that when we find practices and structure that truly serve these students, we do \textbf{much} more than serving these students. We find practices and structures that benefit the wider mathematical community. 

For example, we must confess that the organizing committee had not explicitly thought about the experience of the LGBTQ+ community at ECCO. But our struggles are connected, sometimes in ways that we do not foresee. 
Postdoc Aram Dermenjian \cite{Dermenjian} wrote:
\begin{quote}
\emph{The single-handed biggest reason I loved this conference was the diversity and inclusiveness.
In recent years I felt like the only gay person doing mathematics. 
I've started to feel more and more lonely in my math community. All my friends are amazing and they always try to make me feel welcome, but it's just not the same.\\
Having a community agreement allowed everyone to be open about themselves allowing queer people to be out. 
It gave me the confidence to do something I had never done before. I invited my math friends at ECCO to a gay club. Sure, more than half of us weren't gay, and sure, the gay club wasn't that amazing, but just being out, in a gay club, with mathematicians... was amazing. I felt I belonged.
}
\end{quote}

Professor Viviane Pons \cite{Pons} wrote:
\begin{quote}
\emph{One question arose from the students: why do professors come and teach at ECCO? They saw clearly what was the gain for them, but the reason \textbf{we} would spend time and energy there was not clear to all of them. So let me tell you what I (and the scientific community as a whole) gain from that investment.
I can help shape the academic world to something better and to something I like. Being part of ECCO is one step in this direction, because that is the kind of math community I want.}
\end{quote}

\section{\textsf{The Colombian way}}

Worldwide, people have (accurate or inaccurate) ideas of what is the ``French style," ``Hungarian style," or ``Japanese style" of doing mathematics. In a country that is relatively new to research in mathematics,  perhaps we still have the opportunity to shape what ``Colombian mathematics" might look and feel like. 

Many would argue that mathematics does not distinguish a person's culture or nationality. Unfortunately, this widespread belief has led many of us to feel forced to leave our humanity at the door, and struggle to fit in with the dominant mathematical cultures and practices. 
But our cultures are too rich to be dismissed when we enter mathematics, and dismissing them is a loss to mathematics itself. 

\begin{quote}
\emph{Research has shown that creating learning environments that value and incorporate students', families', and community members' cultural and linguistic strengths into instruction creates a nexus to mathematics cognition. [...] Culture and mathematics learning are intertwined in that they are both transformed through everyday lived experiences and are shaped by those
experiences.}

Michael Orosco and  Naheed Abdulrahim
\cite{OroscoAbdulrahim}
\end{quote}

We are not interested in patriotism, but we \textbf{are} interested in culture and values.
How might we use the cultural practices of Colombian communities to positively influence the cultural practices of Colombian mathematical communities, or at least the cultural practices of ECCO?

Colombians pride ourselves in being good hosts, and making every effort to help our guests feel welcome and comfortable. 

We are proud of our food, our music, our rebusque, and our stories. At ECCO, these all end up playing a central role.

Every ECCO seems to have an improvised fruit tour, where we cross town to visit the local fruit market. It's a bit like going to the zoo, and getting to see all these species that you didn't know existed. One time after dinner, the visiting professors asked students to list all the local fruits they knew;  after about 70 different fruits, everyone got tired of this exercise. 
At the market, Greta bought a huge guan\'abana to share, even though she had no idea what it was. Someone thought it looked like a small alien. Those were not coronavirus times, so no one thought twice about tasting this white, slimy giant, despite the fact that the more people dug in, the blacker it got. 
Eww!!, you might say, and you might be right. But this unplanned experience created a lasting bond. Participants keep connecting around food throughout ECCO; by the end, the foreigners  are showing the locals delicious restaurants and foods that we didn't know.\footnote{For the record, guan\'abanas are heavenly.}

Colombia has a unique salsa culture, where each salsoteca has a wall full of hard-to-find vinyl records from the 60s and 70s -- mostly salsa, cumbia, and West African music -- and this is what the DJ plays all night long. People of all ages dance with their family, their friends, their coworkers, and with any stranger who asks to dance with them. No one does those fancy turns and slimy moves they teach in American salsa lessons. Everyone sings along as they dance.
At ECCO, we organize a visit to the salsoteca with the deepest music collection we can find; C\'esar gives everyone (foreigners and bogotanos alike) dance lessons, the venue gives us maracas and cowbells to play along (and they take them away if our rhythm is not on point), and everyone dances together.
And once the dancing starts, it does not stop. Imagine a conference in which participants work hard during the day, dance during the night and get up fresh and early next morning all over again, for two weeks. Well, dancing drains some energy out of us, but it also fuels us to return to the conference the next day and give it our all. In the final survey, undergraduate student Eliana said: 
\begin{quote}
{\emph{I think that dancing is an important part of ECCO and it changes the whole dynamic in a very positive way. In Colombia mathematics is danced.}}
\end{quote}

Colombians are used to working with a shortage of resources, and we have a strong culture of \textbf{rebusque}: in the face of difficulty, there is always an ingenious solution to be found within our means. 
The first time we organized a summer course in geometric combinatorics for undergraduates, we were advised by foreigners that Colombian students would not have the preparation necessary to understand these topics, that we should teach a basic course in abstract algebra instead. But, for better or for worse, a lack of preparation has never stopped a Colombian from trying to accomplish something. We don't really believe in deficit mindsets. This culture of rebusque probably shapes our conviction that even if you do not have a lot of experience with mathematics, if you are hard-working and resourceful, you can take a class about current research directions from the world experts in the field, learn from it, and contribute to it.

The final activity of ECCO is a panel discussion where we talk about personal issues that most of us struggle with at every stage of our careers, but we rarely or never talk about. Discussion ranges from everyday topics such as  ``What does a typical day in your life look like?'' to more trascendental ones as ``What tools have worked for you to deal with stress, anxiety, or a sense of not belonging in academia?''. We choose a broad range of panelists, from professors to undergrads, from all parts of Colombia and the world.
To be truthful, we are not big fans of math panels in general. Why does this one feel different to us? Perhaps it's that a strong sense of community and trust has been built by the end of ECCO, and this leads to a very honest conversation that does not shy away from strong emotions. Perhaps it's simply that people, and their stories, are really important to us. Daniel wrote:  

\begin{quote}
{\emph{(The panel) is one of the most relevant things at ECCO. Breaking with the idea of math as a selfish and lonely task should be a priority. Talking like human beings, with our emotions and conflicts, is fundamental. Unfortunately this is rarely done in academic events. Congratulations to the organizers for recognizing the need to humanize math and mathematicians. This was a cathartic experience.}}
\end{quote}

This encuentro is an intense mental, physical, and emotional experience. 
We seem to have a tacit agreement to store lots of energy prior to ECCO, and budget a few days of recovery afterwards. We start as a bunch of strangers and end, with tears of joy, promising to keep in touch, planning our next encuentro.

\section{\textsf{What does ECCO wish to be?}} \label{sec:aspirational}

What makes a scholar decide whether to participate in an academic environment? Most academic conferences are conceived as spaces whose main purpose is to advance the creation of knowledge in a very traditional way. This tradition can make such spaces very hostile for many of us given their homogeneity of race, gender, ideology, and experiences.

We love mathematics, and we know that the academic mathematical world will evaluate us and our students by \textbf{what} mathematics we engage with and produce. Thus it is crucial that we offer courses and workshops and are involved in research of a very high quality -- while recognizing that  ``quality" is subjective. We invite professors who are world experts in their fields, excellent lecturers, and compassionate human beings. They give minicourses in some of the most exciting research directions in combinatorics, with strong connections to other fields of mathematics. This strongly influences the kind of research that we and our students end up involved in.

\begin{figure}[h]
\begin{center}
    \includegraphics[width = 8cm]{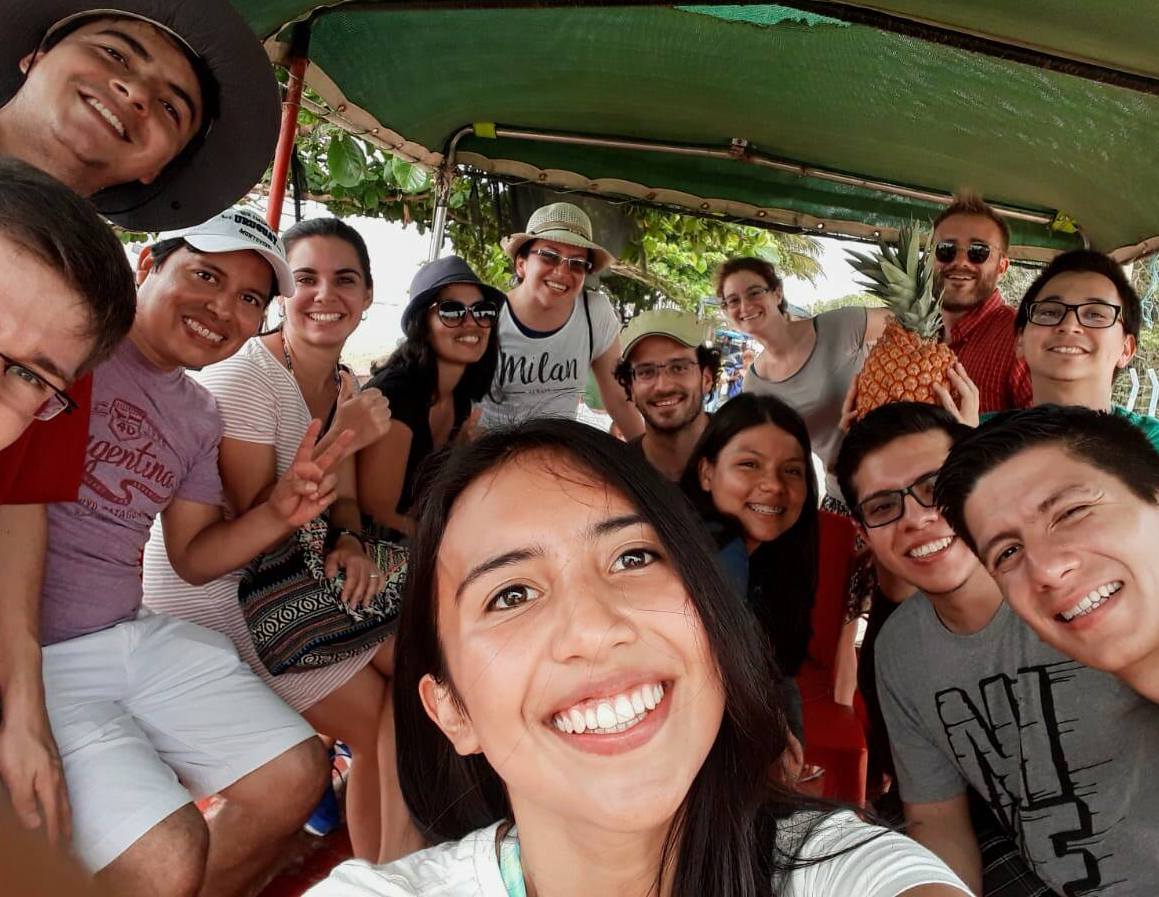}
	\caption{ECCO '18 participants from six countries and five cities in Colombia.}
		\vspace{-.5cm}
\end{center}	
\end{figure}

For us, it is equally crucial to be mindful of \textbf{how} we engage with and produce high quality mathematics. 
Valuing and promoting respect and difference has been essential to the development of ECCO.  We seek to create an environment where each participant is empowered to take the space that belongs to them, and share their voice, ideas, experiences, and world views, inside and outside the classroom. 
More than a conference, ECCO has become a space where learning mathematics is as important as recognizing each other as mathematicians and as individuals. If we are aware of what makes us different mathematically and personally, we can take advantage of these differences to complement each other. 
\begin{quote}
\emph{The difficult, but also essential part, is to value respect and difference positively; not as minor, inevitable nuisances, but as the elements that enrich life and encourage creation and thought.} 

Estanislao Zuleta \cite{Zuleta}
\end{quote}

ECCO brings together a close-knit community of mathematicians who are spread out all over the world. To consolidate this community and institutionalize our efforts, we recently founded an official research group recognized by Colciencias, the Colombian science agency. This has allowed us to systematize and strengthen our research collaborations with each other and with the combinatorics community at large. For example, it has led to the founding of the Seminario Sabanero de Combinatoria\footnote{SeSaCo is a weekly seminar that rotates locations among five universities in Bogot\'a.} in Bogot\'a, and it has strengthened our mathematical connections with ALTENUA\footnote{ALTENUA is a research group in algebra and number theory with a very strong presence in many regions of Colombia.} in Colombia, and with CIMPA\footnote{CIMPA is a nonprofit organization founded in France that promotes research in mathematics in developing countries.} 
worldwide.

Mathematics is not independent from society, and we hope that this mathematical event can be a modest contribution towards the society that we would like to be a part of. Getting more organized institutionally has also helped us initiate, support, and promote various outreach efforts across the country.

\begin{figure}[h]
\begin{center}
    \includegraphics[width = 8cm]{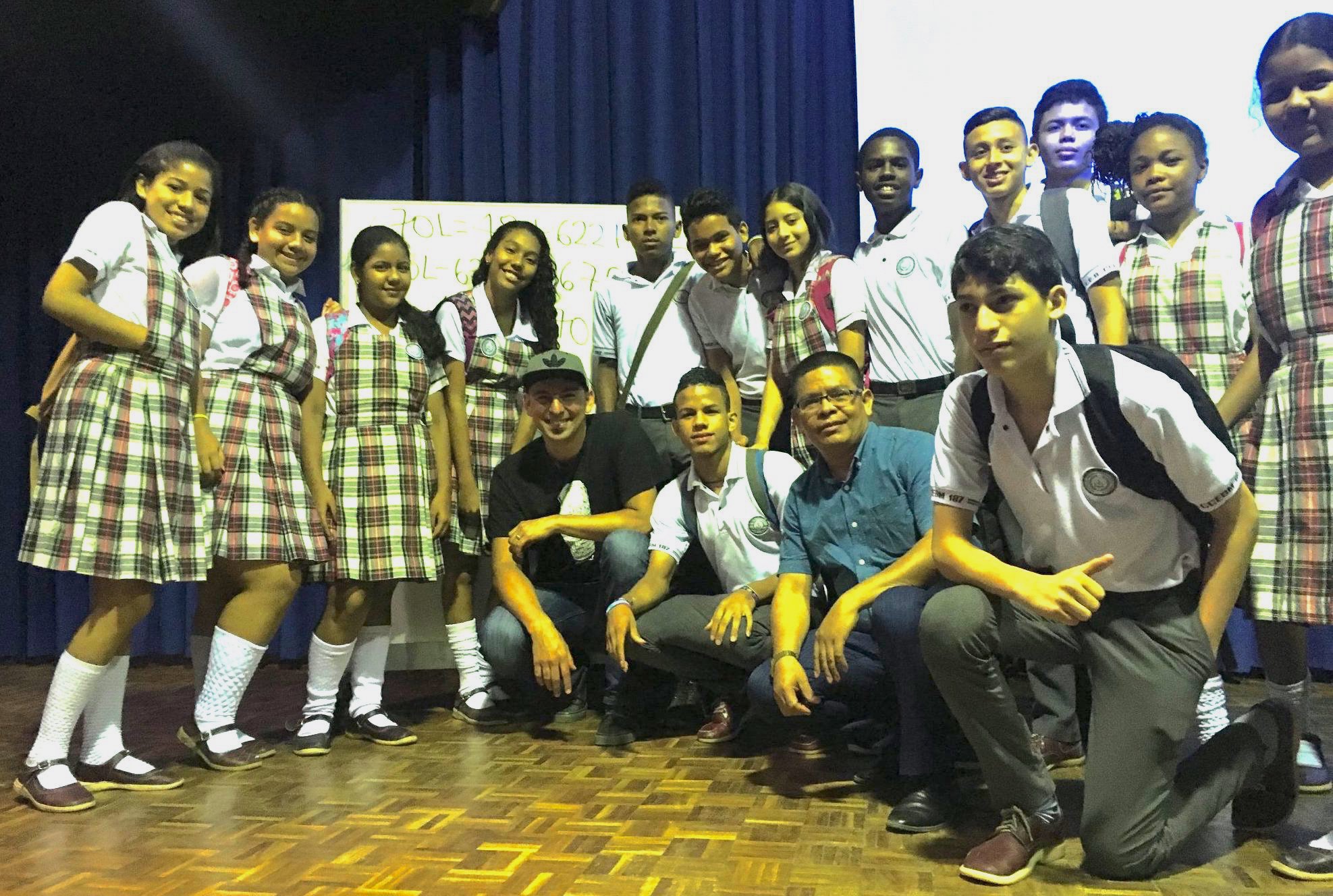}
	\caption{The future: Barranquilla high schoolers.}
	\vspace{-.7cm}
\end{center}	
\end{figure}

During ECCO, we always host a public event in partnership with local initiatives, like Clubes de Ciencia Colombia\footnote{This is a program that helps Colombian researchers abroad lead scientific workshops with students from public high schools in various regions of Colombia. Our work with them is described in \cite{Robots}.}, the science museum Parque Explora in Medell\'{\i}n, and the Nortem\'atica festival for high school students at UniNorte in Barranquilla.

To achieve a more lasting effect, we also challenge ECCO participants to continue to mold their mathematical knowledge, so that it can be used beyond the creation and understanding of science, as a tool of positive societal and human impact.
This has led to the construction of other mathematical communities that celebrate different ways of learning, and help dehomogenize the concept of academia in Colombia. Examples 
include D\'{\i}as de Combinatoria\footnote{D\'{\i}as is a summer school in basic combinatorics, geared towards undergraduates who haven't had access to classes in this area. More than half of D\'{\i}as alumni, from 16 different universities in 9 different cities, went on to attend ECCO.}, and C\'{\i}rculos Matem\'aticos\footnote{C\'{\i}rculos is now a national program that helps high school students from public schools fall in love with mathematics in an inclusive, non-competitive setting.}. 
These are not ``gifted youth programs," which are very valuable but tend to serve the students that already have the most access; our efforts aim to reach the broadest group of participants possible.
Our goal is not to create an army of followers who become mathematicians and perpetuate our idea of academia, but to contribute to a society that appreciates and uses mathematics, and takes advantage of its rich cultural diversity to improve itself. 

\begin{figure}[h]
\begin{center}
    \includegraphics[width = 8cm]{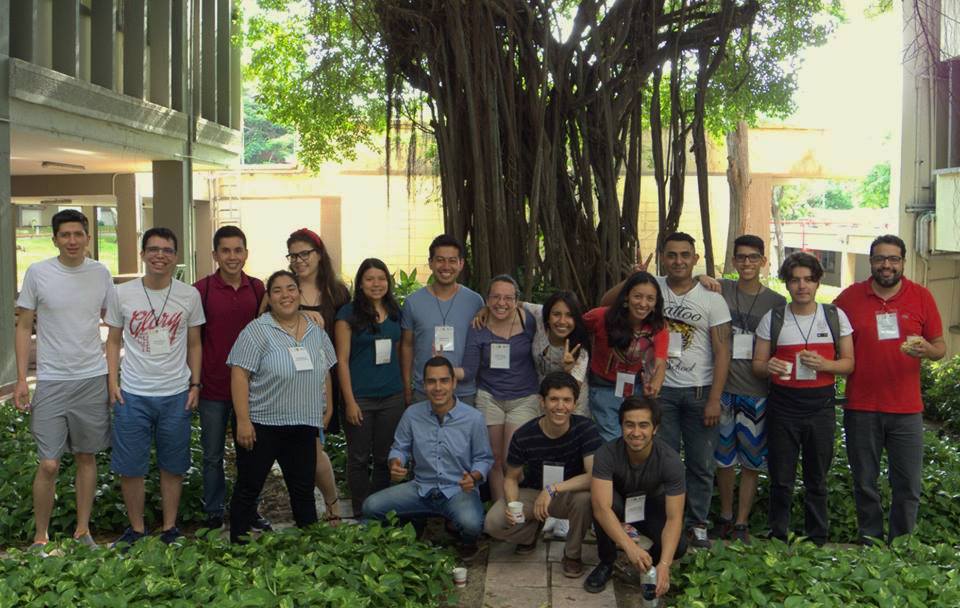}
	\caption{The future: D\'{\i}as '17 alumni at ECCO '18.}
	\vspace{-.5cm}
\end{center}	
\end{figure}

ECCO wishes to help build a strong and dynamic research network that collaborates regionally and internationally, and produces very interesting mathematical work.
We also wish to help create a culture of sharing mathematical knowledge with the public, and using this knowledge to have a positive impact in all sectors of society. 
Finally, we wish to be very mindful of how we do this work, putting our humanity, our values, and the diversity of our cultures at the center of everything that we do. These are the goals that guide our work. In our minds, they are inseparable.

\section{\textsf{Prevention and intervention}}

When we returned home from the most recent ECCO, we received a message from one of our foreign participants. She had very kind words about the conference, particularly about the community agreement and the dancing. However, there were times when the dancing made her uncomfortable, and she felt that some dance partners were not sensitive to her discomfort. 

No conference organizer is happy to receive a message like hers. But it may be a good sign that she felt comfortable enough to share this with us, and trusted that we would listen and take it seriously. We are grateful to her, because her message showed us that we have more work to do to make the community agreement more effective. We have since strengthened  our prevention and intervention protocols: 

1) We have contributed to and incorporated the new protocols of the Commission of Gender and Equity of the Colombian Mathematical Society. \cite{SCM}

2) We have gathered concrete, practical resources for future ECCO organizers on how to respond to potential incident reports, relying on Ashe Dryden's valuable work. \cite{Dryden}

3) We plan on adding a short training on bystander intervention to the schedule; this has proved an effective tool to stop harassment on college campuses. 
Oversimplifying, it proposes: Let's talk about how we treat each other, keep an eye out for each other, and intervene if needed. We need to understand this as a community issue, and not an individual issue, if we want to truly transform the harmful practices that our societies have normalized.

Mathematics has lost too many people -- primarily women and people of color -- to harassment and discrimination, and silence has never protected the victims. Perhaps by sharing with you how we are confronting these problems in our context, we may help you confront them in yours.

\section{\textsf{Keeping the flow}}

\noindent $-$  Fede, why did you start ECCO? Did you have all of this in mind?

\smallskip
\noindent $-$  No, Caro. I can't pretend I was super thoughtful when ECCO started. There was so little activity in combinatorics in Colombia, I just wanted to do \textbf{something}. I did know from the beginning that I wanted it to be really student-centered, collaborative, and unhierarchical. Already in 2003 we had professors and undergraduates working together on problems that neither knew how to solve. 

I was working with students in Colombia and California who were very different from each other, and I knew you could learn a lot from each other. You were also very similar to each other; I thought we could build a space where you all would feel really comfortable and thrive, mathematically and personally. Also, I saw very little understanding and collaboration between Latinx mathematicians from the US and from Latin America, and I thought we could build a bridge.
I knew you and Amanda Ruiz would \textbf{love} working together.

As time went on, more and more young people arrived, and you all brought many great ideas. I loved how it started evolving! 

\smallskip
\noindent $-$  (smirking) ``Organically," like the hipsters say?

\smallskip
\noindent $-$  Organically, like Californians say. We are going to use that word in the article.

\smallskip
\noindent $-$  Noooooo!!!!

\smallskip
\noindent $-$  One thing I did know, early on, is that I could not hold on too tightly; I had to share the decision making, and make myself replaceable as quickly as possible. I have seen many amazing math programs run by the same person for decades; they get stuck! The founders don't want to let go; it's their baby. The participants start respecting the founder too much to point out the serious shortcomings that they see. There is no room for new ideas. It's a real danger. Programs need to evolve.

\smallskip
\noindent $-$  I must confess that I was really worried when you said you weren't going to organize ECCO anymore.

\smallskip
\noindent $-$  Why?

\smallskip
\noindent $-$   Who had the experience, credibility, and recognition you did? Little did I know that we had been learning more than math from you: leadership and community building were growing within each of us. 

\smallskip
\noindent $-$  That's what you think because I'm older. I also feel that I learned this from you all. I think that's what Freire means when he talks about co-constructing knowledge, or experiences: learning from each other as we build it together.

Anyway, you and I know that ECCO has gotten much better since I stopped organizing it. 

\smallskip
\noindent $-$  Hahaha, that's true.

\smallskip
\noindent $-$  What have you tried to bring to ECCO, Caro?

\smallskip
\noindent $-$  I don't know. But I can tell you what I would love to bring to and take from it. ECCO gave me back my hope in academia. I have wanted to run away from this environment many times and the more I get involved in ECCO, the more I realize that it is possible to make use of math knowledge to bring out joy from each of us. That joy then transforms into action which transforms into building more spaces like this. I would love to bring and take that passion that brought me back to wanting to be in academia, at least for the moment. Every time ECCO ends I feel all this energy that wants to be materialized into taking action toward a more inclusive math environment. 

\smallskip
\noindent $-$  I love hearing that. I know that feeling too. But you're being too humble. One of my favorite things about ECCO has been to see you and many others really take ownership of this project, and bring such wonderful new energy, ideas, and perspectives. It has kept ECCO evolving constantly. But more importantly, I believe it has empowered each one of us to take agency and really work towards the mathematical society that we would like to be a part of. 

I am \textbf{so} impressed by your work with Rafa on D\'{\i}as and SeSaCo,
with many others on C\'{\i}rculos 
and in the Comisi\'on, 
by the growing collaborations with ALTENUA and CIMPA,
by the work that Alejo, Ana, Andr\'es, C\'esar, Federico C, Felipe, Jose, Laura, Mandy, Nelly, and \textbf{so} many others are doing with a similar spirit in Colombia and abroad\footnote{Some of these efforts are described in Section \ref{sec:aspirational}.} --  even when I get to march with some of you on the streets when it's needed, I feel so uplifted! I'm not saying ECCO can or should take credit for all of this, but I \textbf{am} saying that ECCO is now one part of a growing ecosystem of interconnected initiatives that I think will have a big impact in the long run. In fact, I think it's the \textbf{only} way we can make a big impact in the long run.
So what do you think we still need to work on? What's next for ECCO? 

\smallskip
\noindent $-$ Ufff... it's funny because I feel that the more we keep doing, the more needs to be done. But don't be surprised if you start seeing Dias de \underline{\hspace{1cm}} schools organized by other people; we have been asked to write a roadmap on how to organize events like this in other areas.

On the other hand, I believe one of the most important things we need to do is to start building trust with regions that have not made part of the things we do. I do not know a definite answer on how that can be done but I am positive that we need to be physically in those places to make that happen. We can start by remembering that context matters, and it is not the same to do science in Bogot\'a or Medell\'{\i}n than to do it in el Pac\'{\i}fico.\footnote{La Regi\'on del Pac\'{\i}fico is a region of great biodiversity and economic importance to Colombia. Its population -- approximately 90\% of which is Afro-Colombian -- has been historically neglected by the national government.}
Context makes people feel heard and recognized. On my end I know I need to educate myself on science from a multicultural perspective. As much as I love all we have been doing with ECCO I feel it is missing one leg. I feel that as long as we are not capable of taking full advantage of the multiculturality of Colombia, we will just continue doing things closer to the traditional way than we think. 

Some tools that may come in handy are the implementation of C\'{\i}rculos and/or D\'{\i}as in other regions. I'd love to see a C\'{\i}rculo Matem\'atico that takes into account the richness and traditions in el Pac\'{\i}fico. I am not saying I want to create it, but rather I would love to contribute to its creation. We could start building bridges with local organizations there to this end.
If  we manage to expand our community to involve members outside of the main cities of Colombia, we will be better off with the future of ECCO. 

ECCO has become, to me, an opportunity to rethink academia and do science my way, and I would love for others to feel the power that has given me. Hopefully in a few years we will stop seeing el Pac\'{\i}fico as that foreign place we are unfamiliar with and only think of when planning holidays, and instead think of it as a place where we draw our role models of academia from.

\begin{figure}[h]
\begin{center}
    \includegraphics[width = 8cm]{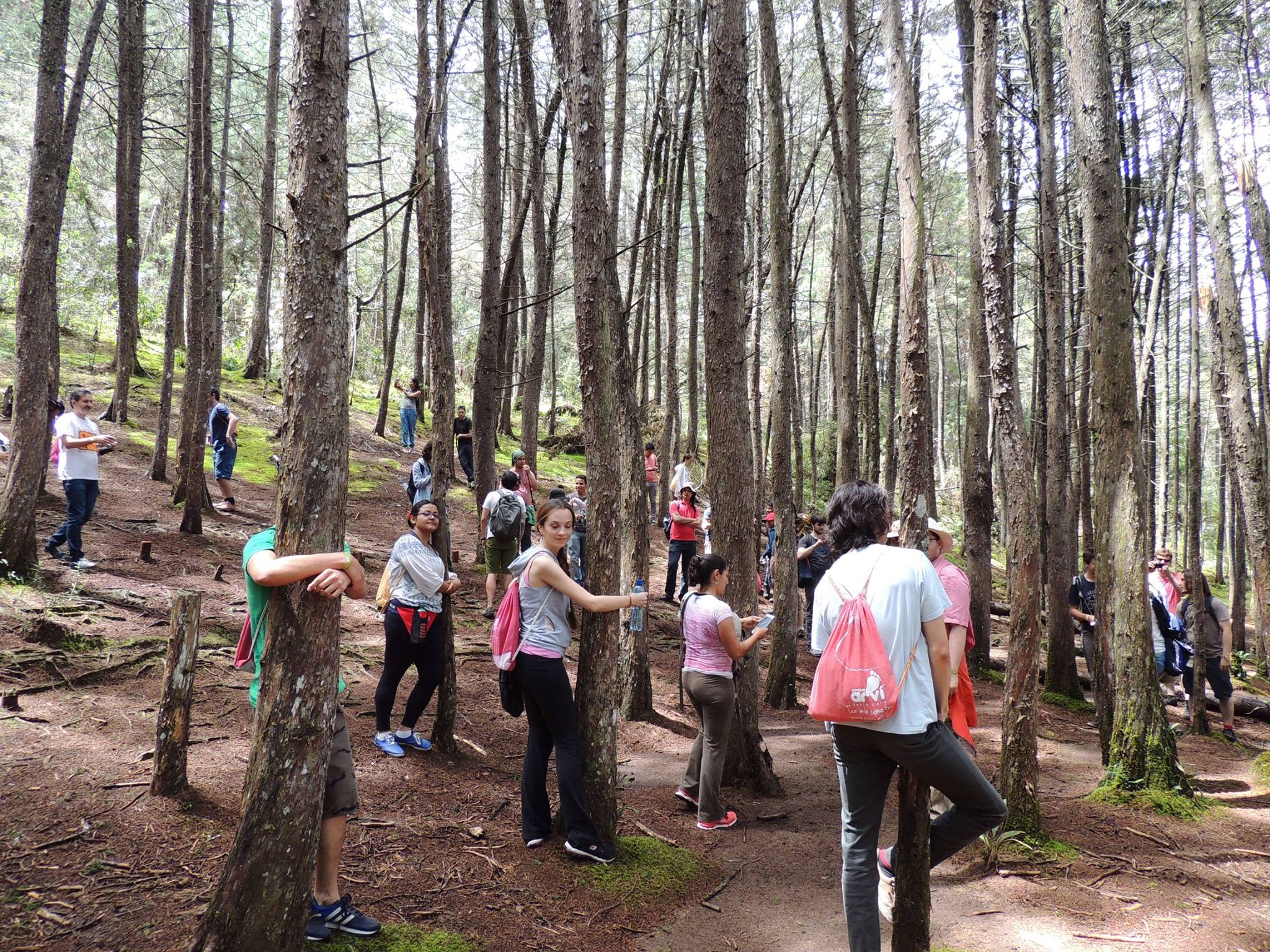}
	\caption{Tree huggers in the Antioquia forest.}
	\vspace{-.7cm}
	\end{center}	
\end{figure}

\section{\textsf{Acknowledgments}}

We would like to extend our warm gratitude to the ECCO family, who have become our mathematical home. They have allowed us to experience true community and belonging in a mathematical space. We also thank the editors of this volume for the invitation to write this article.

FA was supported by NSF grant DMS-1855610, Simons Fellowship 613384, and NIH SF BUILD grant 1UL1MD009608-01. CB was supported by Grant FAPA of the Faculty of Science at the Universidad de Los Andes.

\scriptsize{

}
\end{document}